\documentclass[letterpaper, 10 pt, conference]{ieeeconf}  

\IEEEoverridecommandlockouts                              

\usepackage{graphicx} 
\usepackage{amsmath} 
\usepackage{amssymb}  
\usepackage{amsthm}
\usepackage{algorithm}
\usepackage{algpseudocode}

\newcommand{\R}{\mathbb{R}}
\newcommand{\SOn}{\text{SO}(n)}

\newcommand{\son}{\mathfrak{so}(n)}

\newcommand{\matZ}{\mathcal{Z}}
\newcommand{\bfC}{\mathbf{C}}
\newcommand{\calO}{\mathcal{O}}
\newcommand{\E}{\mathbb{E}}

\newtheorem{theorem}{Theorem}
\newtheorem{assumption}{Assumption}
\newtheorem{lemma}{Lemma}

\title{\LARGE \bf
Tangent Space Parametrization for Stochastic Differential Equations on SO(n)}

\author{Xi Wang$^{1}$ and Victor Solo$^{1}$
\thanks{$^{1}$Xi Wang and Victor Solo is with School of Electrical Engineering \& Telecommunications, UNSW, Sydney, Australia. Victor Solo is the corresponding author. Email: Xi Wang {\tt\small xi.wang14@unsw.edu.au}, Victor Solo {\tt\small v.solo@unsw.edu.au}.  }%
}

\begin{document}

\maketitle
\thispagestyle{empty}
\pagestyle{empty}

\begin{abstract}
In this paper, we study the numerical simulation of stochastic differential equations (SDEs) on the special orthogonal Lie group $\SOn$. We propose a geometry-preserving numerical scheme based on the stochastic tangent space parametrization (S-TaSP) method for state-dependent multiplicative SDEs on $\SOn$. The convergence analysis of the S-TaSP scheme establishes a strong convergence order of $\calO(\delta^{\frac{1-\epsilon}{2}})$, which matches the convergence order of the previous stochastic Lie Euler-Maruyama scheme while avoiding the computational cost of the exponential map. Numerical simulation illustrates the theoretical results.
\end{abstract}

\section{INTRODUCTION}
Over the past few decades, stochastic differential equations (SDEs) evolving on Lie groups have played an important role in mathematics~\cite{hsu2002stochastic}. Among them, SDEs on the special orthogonal group $\SOn$ have attracted significant research attention due to their wide range of applications in robotics,  molecular dynamics, and quantum control~\cite{chirikjian2011stochastic,xiang2020simultaneous}.

For existing applications, it is important to develop geometry-preserving numerical schemes (GPNS) for SDEs on $\SOn$. The purpose of a GPNS is to generate approximate solutions that remain entirely within the Lie group. Therefore, numerical schemes for Euclidean SDEs, such as those in \cite{platen1999introduction} and \cite{higham2001algorithmic}, cannot be directly applied to SDEs on Lie groups, as Euclidean methods do not account for the geometric constraint, 
and lead to trajectories drifting away from the manifold
over long time scales.

%
It is natural to try to extend  ordinary differential equation (ODE) numerical methods
on $\SOn$~\cite{hairer2006geometric} to SDEs on $\SOn$. This leads to numerical analysis schemes based on Stratonovich-SDEs (S-SDEs). But recently ~\cite{solo2024stratonovich} showed that when the diffusion term is state dependent the Stratonovich scheme converges to the wrong SDE. The problem is that numerical schemes converge to Ito-SDEs (I-SDEs) and so numerical schemes must be based on the I-SDE not the S-SDE. Another way to see this is that ODEs evolve in the tangent plane whereas I-SDEs also contain a `pinning-drift'~\cite{solo2024stratonovich}  in the normal to the Lie-Group.


Regarding GPNS for SDEs on $\SOn$, \cite{faou2009conservative} and \cite{armstrong2022curved} introduce GPNS for SDEs on general manifolds, where their approaches rely on either projection onto the manifold or solving Euclidean ODEs. As a result, they lack closed-form solutions.  Specifically for $\SOn$, \cite{piggott2016geometric} proposed a stochastic Lie EM (SL-EM) scheme with a strong convergence order of $\mathcal{O}(\delta^{\frac{1}{2}})$ when the diffusion term is constant. \cite{muniz2022higher} developed a stochastic Lie RKMK method for state-independent SDEs on $\SOn$, and proved the strong convergence order $\mathcal{O}(\delta^k)$ for the $k$-th order Runge–Kutta method. Recently, \cite{solo2024stratonovich} established a strong convergence rate of $\mathcal{O}(\delta^{\frac{1-\epsilon}{2}})$ for state-dependent SDEs on $\mathrm{SO}(n)$, where $\epsilon$ can be arbitrarily small.

All of the above works rely on the Magnus expansion of the exponential map on Lie groups~\cite{marjanovic2018numerical}. However, explicitly computing the exponential map, such as using the Rodrigues formula on $\SOn$, can be exponentially ill-conditioned~\cite{gallier2003computing}. Moreover, in practice, the exponential map is often approximated using the Padé approximation with scaling and squaring~\cite{piggott2016geometric}, which involves computing a rational function of a matrix. Therefore, computing the exponential map on Lie groups can be computationally expensive.

An alternative approach avoiding computing the exponential map is the tangent space parametrization (TaSP), which originally emerged as a method for solving ODEs on Riemannian manifolds. Convergence results for TaSP applied to geometric SDEs remain quite limited. \cite{wang2021numerical} proposed a TaSP scheme for SDEs on Stiefel manifolds, but no convergence results were provided. \cite{solo2021convergence} showed convergence 
when the Stiefel manifold reduces to a sphere.

Motivated by this research gap, we propose the TaSP method as a numerical scheme for SDEs on $\SOn$ and analyze its convergence rate. The contributions of this paper are summarized as follows:
\begin{itemize}
    \item We introduce the stochastic tangent space parametrization (S-TaSP) for state-dependent multiplicative SDEs on $\SOn$. In contrast to the S-TaSP method on the Stiefel manifold, where the iteration is based on numerically solving quadratic matrix Riccati equations. The Lie group structure of $\SOn$ provides a simple and explicit framework to iterate.
    \item We provide a convergence analysis of S-TaSP on $\SOn$, showing a strong covergence rate of $\calO(\delta^{\frac{1-\epsilon}{2}})$, which matches the convergence order of previous work while avoiding the computation of the exponential map.  
    \item We validate our theoretical results through numerical simulation, illustrating the effectiveness of our method.  
\end{itemize}
The remainder of the paper is organized as follows. In Section II, we reviews the geometry of $\SOn$ and the formulation of SDEs on $\SOn$. In Section III, we introduces the proposed S-TaSP scheme for $\SOn$. In Section IV, we presents the convergence analysis of the S-TaSP scheme. Section V provides numerical simulations to validate the theoretical results, and Section VI concludes the paper. Proof details are included in the Appendix.

\section{PROBLEM STATEMENT}
In this section, we first review the geometry of the special orthogonal group $\SOn$. We then present the multiplicative SDE on $\SOn$ and formally introduce our problem.

\subsection{Geometry of $\SOn$}
The special orthogonal group $  \SOn = \{ R \in \R^{n \times n} \mid R^T R = I, \det(R) = 1 \}$ is the set of all $n \times n$ orthogonal matrices with determinant $1$. The associated Lie algebra $\son = \{ \matZ \in \R^{n \times n} \mid \matZ^T = -\matZ \}$ is defined as all $n\times n $ skew-symmetry matrix.


Given a point $R \in \SOn$, we can write the tangent space at $R$ as $T_R \SOn = \{ R \matZ \mid \matZ \in \son \}$. The normal space of $\SOn$ at $R$ is defined as the orthogonal complement of  $T_R \SOn$ in $\R^{n \times n}$, which can be written as $N_R \SOn = \{ R \bfC \mid \bfC \in \R^{n \times n}, \bfC = \bfC^T\}$. Projecting a matrix $A \in \R^{n \times n}$ onto the tangent space $T_R \SOn$ can be done by the projection operator $\Pi_R(X) = \frac{1}{2} (A - R^T A^T R)$.

\subsection{Multiplicative SDE on $\SOn$}
The multiplicative S-SDE on $\SOn$ is given by
\begin{align}\label{eq:SDE-S}
    dR = RB_{o,S}(R,t)dt + \sum_{j = 1}^d R B_{j}(R,t)\circ dW_j,
\end{align}
where $B_{o,S} , B_j$ are $n \times n$ matrices, $W_j$ are independent Brownian motions with variance $dt$ and $\circ$ denotes the Stratonovich integral. To ensure that the solution $R(t)$ remains in $\SOn$, it is necessary and sufficient that $B_{o,S}(R,t)$ and all $B_j(R,t)$ are skew-symmetric. Using the Stratonovich-to-It\^{o} conversion formula~\cite{chirikjian2011stochastic,solo2024stratonovich}, the equivalent I-SDE in~\eqref{eq:SDE-S} is given by
\begin{align}\label{eq:SDE-I}
     dR = RB_{o,I}(R)dt + \sum_{j = 1}^d R B_{j}(R) dW_j
\end{align}
where the It\^{o} drift term $B_{o,I}$ is expressed as
\begin{align*}
\begin{cases}
     B_{o,I} = B_{o,S} + \frac{1}{2} \sum_{j=1}^d K_j + \frac{1}{2} \sum_{j = 1}^d B_j^2, \\
     K_j = \sum_{r,s = 1}^n \frac{\partial B_j}{R_{rs}} (R B_j)_{rs}.
\end{cases}
\end{align*}

From~\eqref{eq:SDE-I}, note that all $ K_j $ are skew-symmetric, while the term $ \frac{1}{2} \sum_{j = 1}^d B_j^2 $ is symmetric. Therefore, the It\^{o} drift term $ R B_{o,I} $ contains normal components. This normal component highlights a fundamental difference between SDEs and ODEs on $\SOn$, where in the latter, $ K_j = B_j = 0 $. Additionally, it introduces challenges in developing numerical schemes of SDEs on $\SOn$.  

\subsection{Problem Definition} In this paper, we aim to develop geometry-preserving numerical schemes (GPNS) for SDE~\eqref{eq:SDE-I} on $\SOn$ based on tangent space parametrization. A GPNS involves discretizing a time interval $[0, T]$ into $M$ equal time steps, each of size $\delta = T / M$. At each time step, the algorithm updates the state variable $\hat{R}_t$ based on the drift and diffusion terms of the SDE, yielding an approximate solution $\hat{R}(t): [0, T] \to \SOn$. A GPNS for the SDE in~\eqref{eq:SDE-I} is said to achieve a strong convergence order $p$ if the expected error satisfies

\[
    \mathbb{E} \left[ \max_{0 \leq t \leq T} \| R(t) - \hat{R}(t) \|^2 \right]^{\frac{1}{2}} \leq C \delta^p,
\]

where $C$ is a constant independent of the step size $\delta$. 

\section{STOCHASTIC TaSP NUMERICAL SCHEME}
In this section, we propose the S-TasP method for the SDE~\eqref{eq:SDE-I} on $\SOn$. Suppose we are at the point $R_0 \in \SOn$ and aim to move in the direction of the tangent vector $R_0 \matZ$, where $\matZ \in \son$. A naive update of the form $ R_0 + R_0 \matZ $ will generally not preserve the constraint for $ \SOn $. The key idea of the TaSP method is to introduce a normal correction term $ R_0 \bfC \in N_{R_0} \SOn $, such that the updated point $ R = R_0 + R_0 \matZ + R_0 \bfC $ remains on the Lie group $ \SOn $.

In~\cite{wang2021numerical}, the authors studied the SDEs on Stiefel manifold, where the correction $ \bfC $ is obtained by solving a matrix Riccati equation numerically. In this work, we show that the Lie group structure allows a simple and explicit correction term on $ \SOn $, as stated in the following theorem.
\begin{theorem}\label{thm:correction}
    Let $ R_0 \in \SOn $ and $ \matZ \in \son $. If $ I - \matZ^T \matZ $ is positive definite, then the normal correction term $ \bfC_0 \in N_{R_0} \SOn $ that ensures $ R = R_0 + R_0 \matZ + R_0 \bfC \in \SOn $ is given by  
    \[
    \bfC = \sqrt{I - \matZ^T \matZ} - I,
    \]
    where $ \sqrt{A} $ denotes the unique positive definite square root of the positive definite matrix $ A $.
\end{theorem}

The correction term proposed in Theorem~\ref{thm:correction} is a key component for establishing the S-TasP method for SDEs on $\SOn$. For a small time interval $[0, \delta]$ and initial state $R_0$, the real solution of~\eqref{eq:SDE-I} can be parametrized by the tangent vector $\matZ(t)$ as
\begin{align}\label{eq:real-sol}
    R(t) & = R_0(I + \matZ(t) + \bfC(t)), \forall t \in [0, \delta].
\end{align}
Taking differentials of~\eqref{eq:real-sol} we obtain
\begin{align}\label{eq:diff}
    dR = R_0 d\matZ + R_0 d\bfC 
\end{align}
Projecting~\eqref{eq:diff} onto the tangent space $T_{R_0} \SOn$, we have
\begin{align}\label{eq:proj}
    \Pi_{R_0}(dR) = \Pi_{R_0}(R_0 d\matZ + R_0 d\bfC) = R_0 d\matZ.
\end{align}
The last equality holds since the $\bfC$ is symmetric and hence $R_0d\bfC \in N_{R_0} \SOn$, while $R_0d\matZ \in T_{R_0} \SOn$.
Substituting~the SDE~\eqref{eq:SDE-I} into \eqref{eq:proj}, we obtain
\begin{align}\label{eq:TaSP}
    d \matZ = & R_0^T \Pi_{R_0}(RB_{o,I} + \sum_{j = 1}^d R B_j dW_j)  \nonumber  \\
            =  & \left(B_{o,I} - \frac{1}{2} \sum_{j=1}^d B_j^2\right) dt + \sum_{j = 1}^d B_j dW_j.
\end{align}
which is a standard Euclidean SDE. In this case, 
the S-TaSP scheme applies the one-step Euler 
method at time $t = 0$ for the SDE~\eqref{eq:TaSP}, yielding
\begin{align*}
    \matZ_0 = \left(B_{o,I} - \frac{1}{2} \sum_{j=1}^d B_j^2\right) \delta + \sum_{j = 1}^d B_j \sqrt{\delta} \epsilon_j,  \epsilon_j \sim \mathcal{N}(0,1).
\end{align*}
Then the S-TaSP scheme updates the state variable $R_1$ by
\begin{align*}
    R_1 = R_0(I + \matZ_0 + \bfC_0) = R_0(\matZ_0 + \sqrt{I - \matZ_0^T \matZ_0})
\end{align*}
after which we restart the process from the new point $R_1$. The whole process of TaSP can be summarized in Algorithm~\ref{alg:TasP}.

\begin{algorithm}[H]
    \caption{S-TaSP for Multiplicative SDEs on $\SOn$}
    \begin{algorithmic}[1]
    \State \textbf{Input:} Length $\delta$, final time $T$, initial state $R_0 \in \SOn$.
    \State Divide the interval $[0, T]$ into $M$ intervals with length $\delta$.
    \For{$m = 0, \dots, M$}
        \State Draw $\epsilon_{j,m} \sim \mathcal{N}(0, 1)$ and compute
        \begin{align*}
            &\matZ_m=  \big( B_{o,I} - \frac{1}{2} \sum_{j=1}^d B_j^2 \big) \delta + \sum_{j=1}^d B_j \sqrt{\delta} \epsilon_{j,m}. 
        \end{align*}
        \If{$I - \matZ_m^T \matZ_m$ is not positive definite}
            \State Repeat Step 4.
        \Else
            \State Set the correction term $\bfC_m =  \sqrt{I- \matZ_m^T \matZ_m}-I$.
            \State Update $\hat R_{m+1} = \hat R_m(I + \matZ_m + \bfC_m) $
        \EndIf
    \EndFor
    \State \textbf{Output:} Trajectory $\hat R(t) = \hat R_m, (m-1)\delta < t \le m\delta$.
    \end{algorithmic}
    \label{alg:TasP}
\end{algorithm}

\section{CONVERGENCE ANALYSIS}
To analyze the convergence of the S-TaSP scheme, we impose the following assumptions, which are adapted from standard assumptions for SDE numerical scheme on Euclidean spaces or Lie groups~\cite{platen1999introduction,solo2024stratonovich}.
    \begin{assumption}\label{asm:Lip}
        The term $B_{o,I}$ and $B_j, j \in [d]$ are $L$-Lipschitz continuous on $\SOn$, i.e., there exist an $L\in \R$ such that
        \begin{align*}
        \begin{cases}
             B_{o,I}(R) - B_{o,I}(R') \le L\|R-R'\|_F, \\
            B_{J}(R) - B_{J}(R') \le L\|R-R'\|_F, j \in [d].
        \end{cases}
        \end{align*}
    \end{assumption}
    \begin{assumption}\label{asm:bound}
        The term $B_{o,I}$ and $B_j, j \in [d]$ are bounded on $\SOn$ with bounds $C$.
    \end{assumption}

Then we carry out our convergence analysis by introducing the immediate sequence $ R_{\delta}$ as
\begin{align*}
    R_{\delta}(t) - R_0 = & \int_0^{t} \hat R(t)B_{o,I}(\hat R(t)) dt \\
    & + \int_0^{t} \hat R(t) \sum_{j=1}^d B_j(\hat R(t)) dW_j(t)
\end{align*}
In this case, the convergence analysis of the S-TaSP scheme is divided into two parts
\begin{align*}
    & \E \left[ \sup_{0 \le t \le T}\| \hat R(t) - R_{\delta}(t) \|_F^2 \right]  \\
    & \le 2 \E \left[ \sup_{0 \le t \le T}\| a_{\delta}(t) \|_F^2 \right] + 2 \E \left[ \sup_{0 \le t \le T}\| \nu_{\delta}(t) \|_F^2 \right], 
\end{align*}
where
\begin{align*}
    a_{\delta}(t) = R_{\delta}(t) - \hat R(t), \quad \nu_{\delta}(t) = \hat R(t) - R(t).
\end{align*}
We analyze the above two terms in the following lemmas.
\begin{lemma}\label{lem:1}
    Under Assumptions~\ref{asm:Lip} and~\ref{asm:bound}, the error between $\hat R(t)$ and the $R_{\delta}(t)$ satisfies
    \begin{align*}
        \E \left[ \sup_{0 \le t \le T}\| \hat R(t) - R_{\delta}(t) \|_F^2 \right] \le \calO(\delta^{1-\epsilon}),
    \end{align*}
    for any arbitrary small $\epsilon > 0$.
\end{lemma}
The proof of Lemma~\ref{lem:1} is given in Appendix. 
\begin{lemma}
    Under Assumptions~\ref{asm:Lip} and~\ref{asm:bound}, we have
    \begin{align*}
        \E \left[ \sup_{0 \le t \le T}\| R_{\delta}(t) - R(t) \|_F^2 \right] \le \calO(\delta^{1-\epsilon}).
    \end{align*}
\end{lemma}
\noindent \textit{Proof:} It follows Result IIIb in~\cite{solo2021convergence}. \hfill $\blacksquare$ \\
Combining the two lemmas, we have the main result of the convergence order $\calO(\delta^{\frac{1-\epsilon}{2}})$ of the S-TaSP scheme as follows.
\begin{theorem}
    Under Assumptions 1 and 2, the S-TaSP scheme converges to the real solution $R(t)$ of the SDE~\eqref{eq:SDE-I} on $\SOn$ with strong order $\frac{1-\epsilon}{2}$, i.e.,
    \begin{align*}
        \E \left[ \sup_{0 \le t \le T}\| \hat R(t) - R(t) \|_F^2 \right] \le \calO(\delta^{1-\epsilon}),
    \end{align*}
    for any arbitrary small $\epsilon > 0$.
\end{theorem}

\section{Numerical simulation}
In this section, we provide numerical simulation to 
illustrate the accuracy of the S-TaSP scheme for solving 
the multiplicative SDE on $\SOn$ and compare the 
S-TaSP scheme with the stochastic Lie EM (SL-EM) 
method based on Magnus expansion in~\cite{solo2024stratonovich}.

\subsection{Brownian Motion on $\SOn$}
\begin{figure}[tpb]
    \centering  
    \vspace{1em}
   \includegraphics[width=0.4\textwidth]{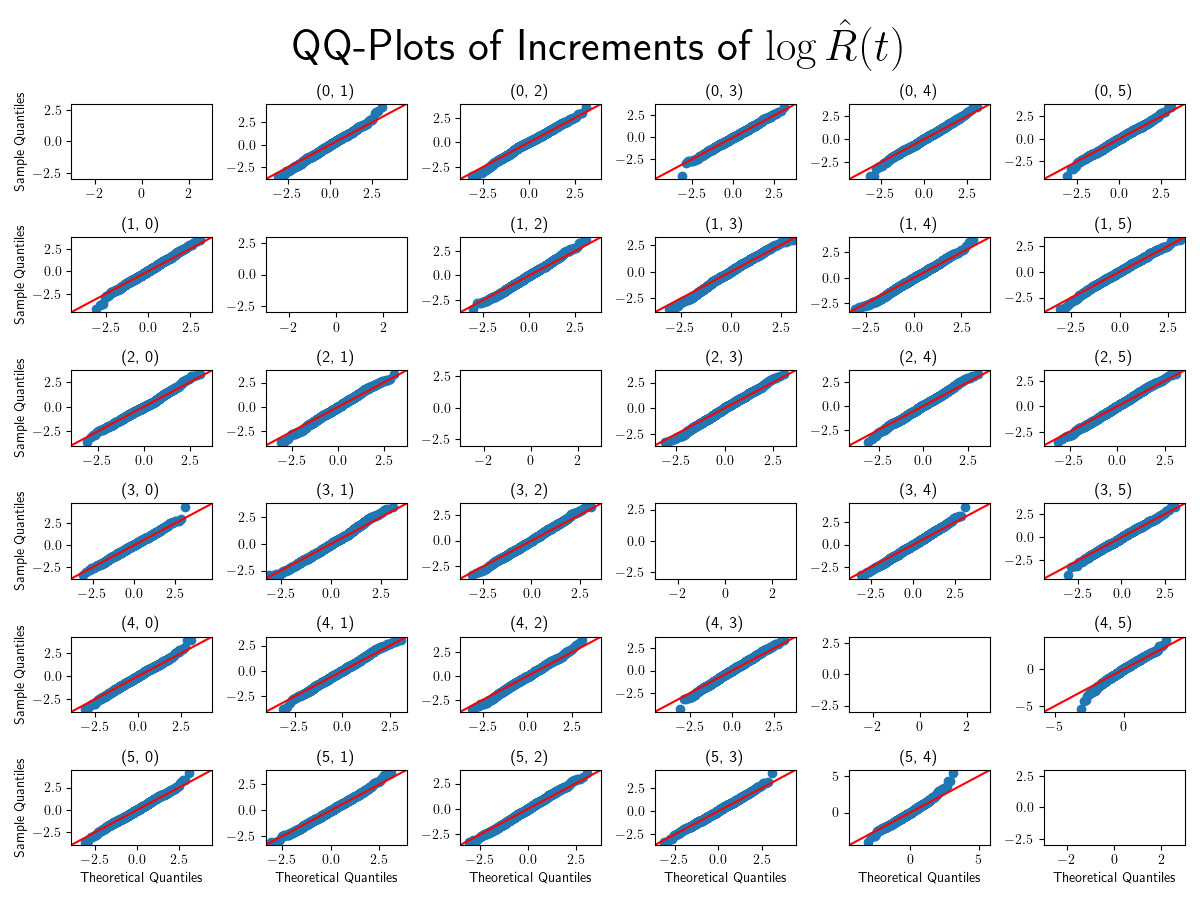}
    \caption{Normality test of $\log \hat R(t)$}
    \label{fig:bm_normal}
\end{figure}

\begin{figure}[t]
    \centering
        \vspace{1em}
    \includegraphics[width=0.4\textwidth]{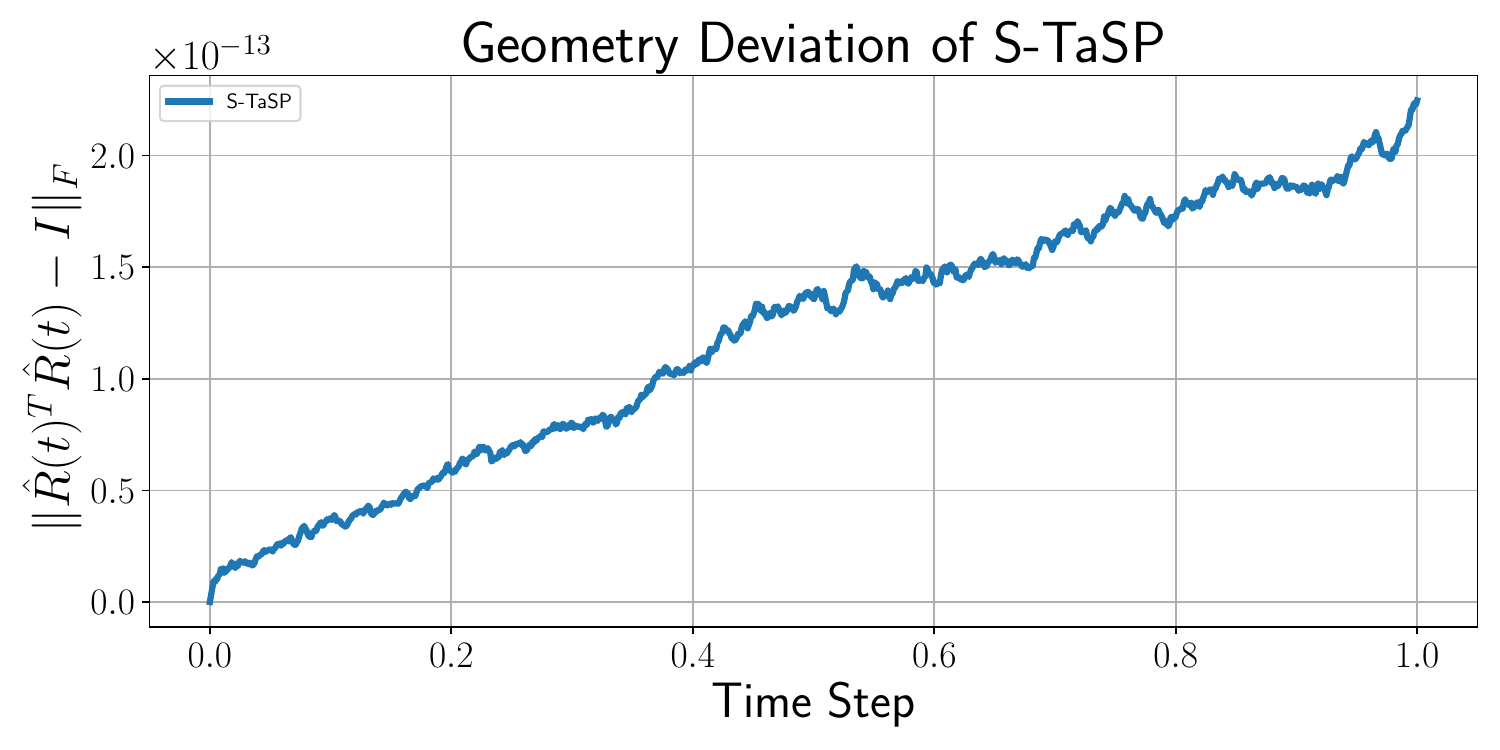}
    \caption{Geometry Preservation of S-TaSP for Brownian Motion}
    \label{fig:bm_geo_diff}
\end{figure}
We consider the following SDE on $\SOn$ 
\begin{align}\label{eq:son_bm}
    dR = -\frac{n-1}{4} RI \, dt + R \sum_{j=1}^{\frac{n(n-1)}{2}} \frac{1}{\sqrt{2}} E_{j} \, dW_j,
\end{align}  
where $E_j$ represents the $j$-th canonical basis vector in $\mathfrak{so}(n)$. In this setup, $B_{o,I}(R) = -\frac{n-1}{4} I$ and $B_j(R) = \frac{1}{\sqrt{2}} E_j$. 

According to \cite{piggott2016geometric}, setting $R(t) = R_0 \exp(\Omega(t))$, we have $d\Omega = d\exp^{-1}(V_t) \approx V_t$, where $V_t = \sum_{j=1}^{\frac{n(n-1)}{2}} E_{j} \frac{1}{\sqrt{2}} dW_j $, which indicates that $\log (R(0)^T R(t))$ is an $n \times n$ skew-symmetric matrix, and every off-diagonal element is an approximate Brownian motion with variance $\frac{1}{2}dt$ in a small time scale.

We tested \eqref{eq:son_bm} using the TaSP scheme for $n = 6$ over the time interval $[0,1]$. In this experiment, we initialized $R_0 = I$ and used a time step $\delta = 0.001$. Figure~\ref{fig:bm_normal} displays the QQ-plot of the elements of $\sqrt{\frac{\delta}{2}} \log (\hat R(t))$ against the standard normal distribution. The results illustrate that the increments of $\sqrt{\frac{\delta}{2}} \log (\hat R(t))$  are normally distributed, which illustrates the correctness of the S-TaSP scheme on $\SOn$.
    


Besides, we investigate the geometric preservation of the S-TaSP. Figure~\ref{fig:bm_geo_diff} plots $\|R^T R-I\|_F$ versus the time $t$. The result in Figure~\ref{fig:bm_geo_diff} shows the geometric deviation of $R(t)$ is maintained at an $O(10^{-13})$ level, which decipts the our S-TaSP scheme preserve the geometric feature. 

\subsection{State-Dependent Noise}
We consider the following state-dependent Stratonovich SDE on $\SOn$:  
\begin{align}\label{eq:sgd-s}
    dR =  -\nabla F(R)  dt +  \tau(R) \sum_{j=1}^{\frac{n(n-1)}{2}} R E_j \circ dW_j, 
\end{align}  
where $F(R):\SOn \to \R$ is a function, and $E_j$ also denotes the $j$-th canonical basis vector in $\mathfrak{so}(n)$. The equation~\eqref{eq:sgd-s} represents the stochastic gradient descent of $F$ on $\SOn$ with random noise.

\begin{figure}[t]
    \centering
    \includegraphics[width=0.4\textwidth]{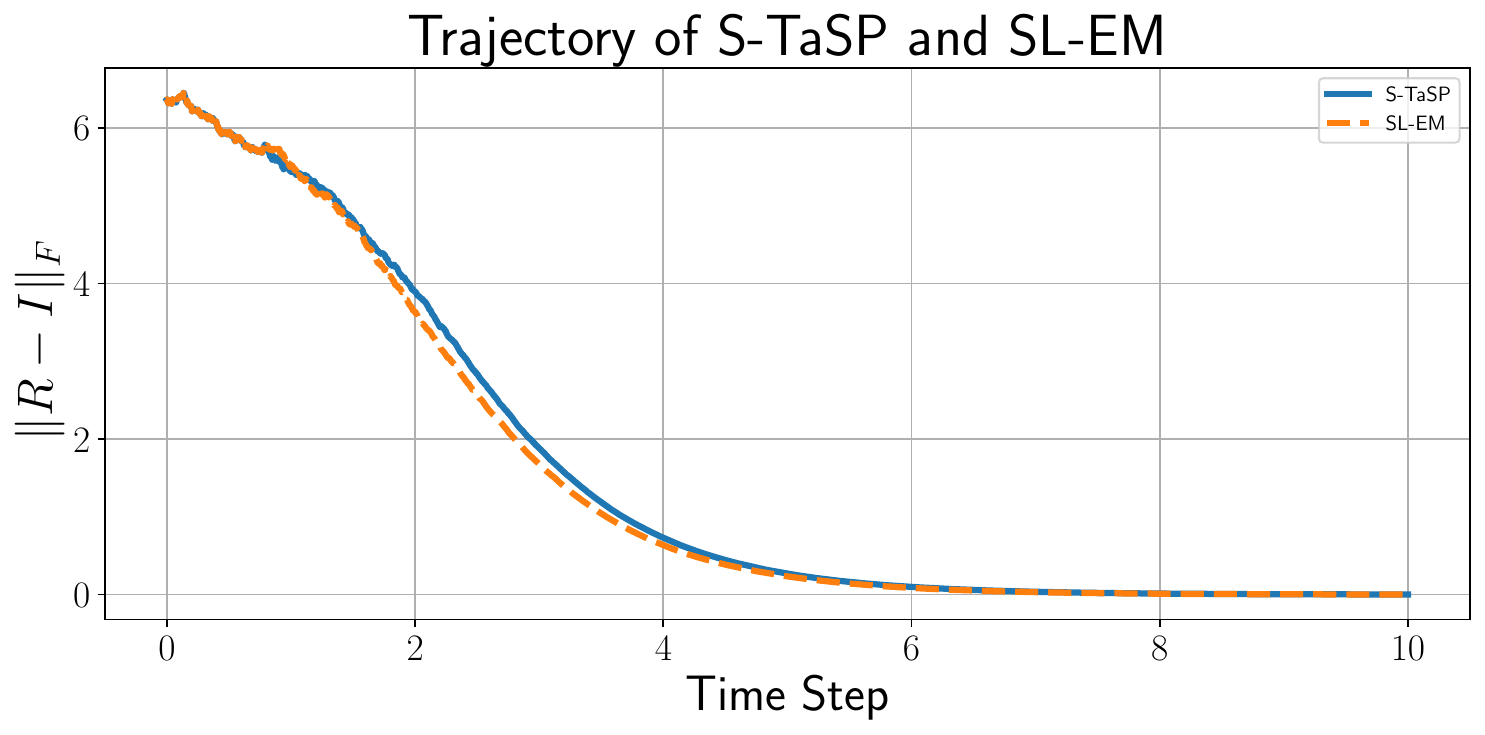}
    \caption{Trajectory Comparison between S-TaSP and SL-EM }
    \label{fig:omega_traj}
\end{figure}

\begin{figure}[t]
    \centering
        \vspace{1em}
    \includegraphics[width=0.4\textwidth]{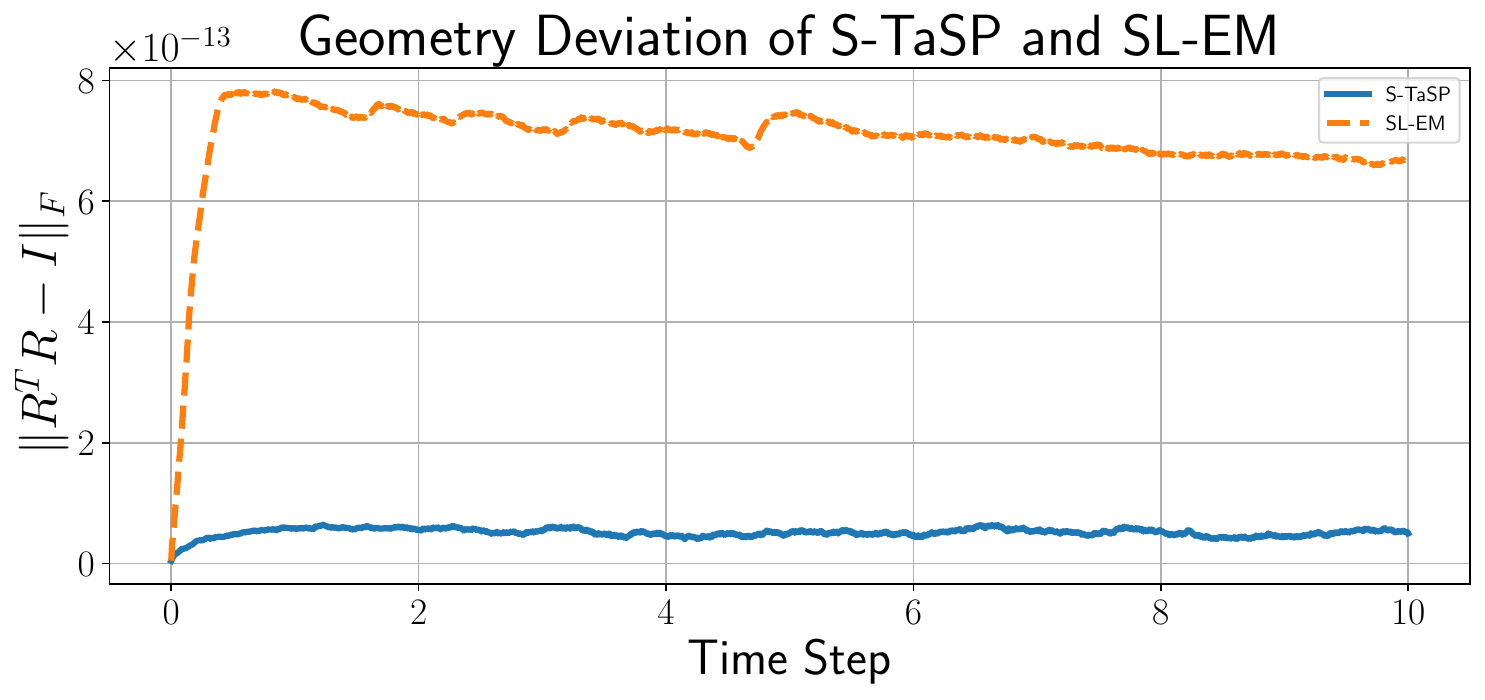}
    \caption{Geometry Preservation of S-TaSP and SL-EM}
    \label{fig:sdn_traj}
\end{figure}

\begin{figure}[t]
    \centering
    \includegraphics[width=0.4\textwidth]{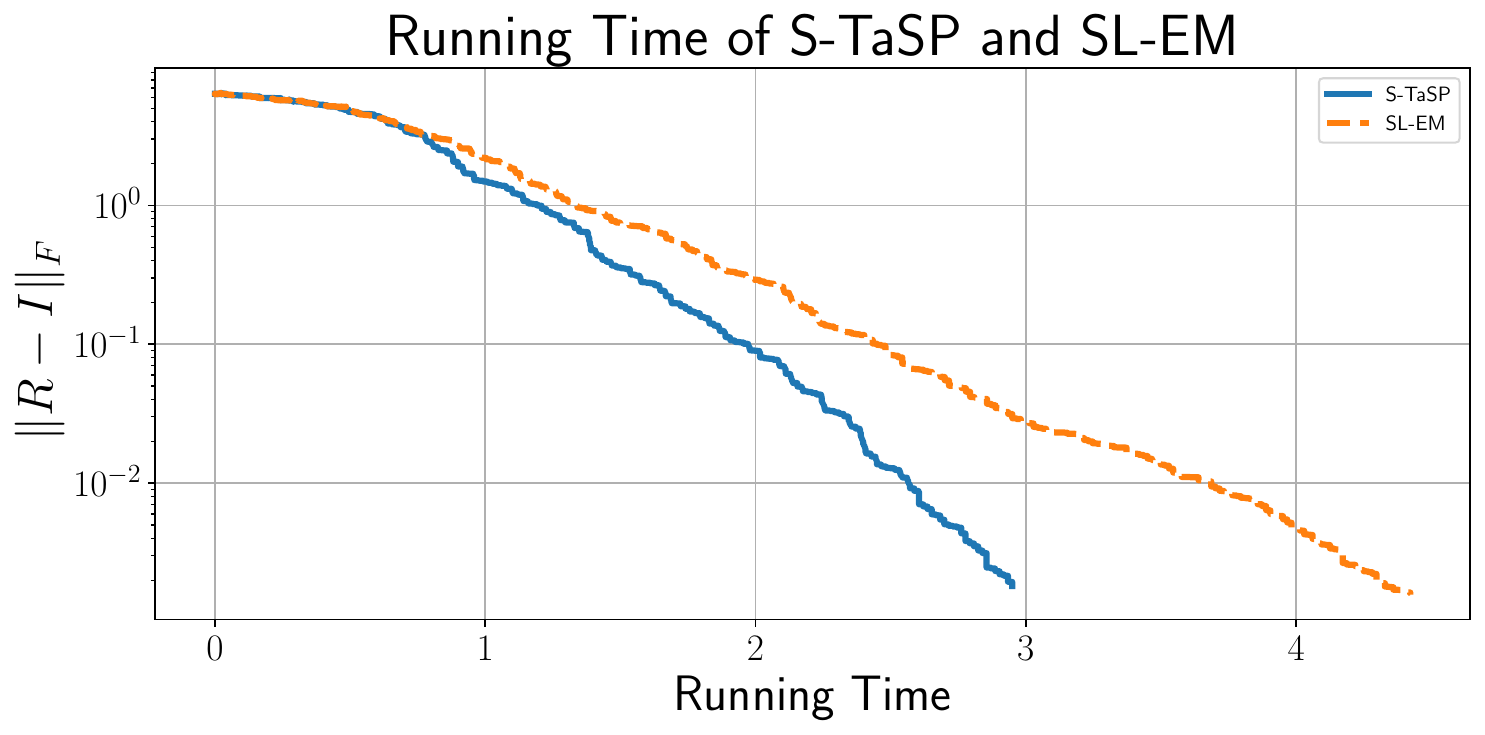}
    \caption{Running Time Comparison between S-TaSP and SL-EM}
    \label{fig:time}
\end{figure}
In our experiment, we set $F(R) =\frac{1}{2} d_{\SOn}^2(R, I)$, where $d_{\SOn}(R_1, R_2) := \|\log(R_1^T R_2)\|$ denotes the Riemannian distance between $R_1$ and $R_2$. The gradient is given by $\nabla F(R) = R\log(R)$. Additionally, we define $\tau(R) = \frac{1}{2} - \frac{1}{2n} \text{tr}(R)$, which represents the noise term. This noise gradually vanishes as $R$ approaches the optimal point $I$. Therefore, the trajectory of~\eqref{eq:sgd-s} will eventually converge to $ I $. In the cases, the corresponding It\^{o} SDE can be written as 
\begin{align*} 
    dR = RB_{o,I}(R)dt + \sum_{j = 1}^d R B_{j}(R) dW_j,   
\end{align*}
where $B_{o,I}= -\log(R) + \frac{\tau(R)}{2n} (R-R^T) - \tau^2(R)\frac{n-1}{2}I$ and $B_j = \tau(R) E_j $. 

We conduct the simulation with $n = 50$ using the S-TaSP and SL-EM schemes over the time interval $[0,10]$ with a step size of $\delta = 0.001$. The square root $\sqrt{I-\matZ^T\matZ}$ is approximated by its $5$-th Taylor expansion. Figure~\ref{fig:omega_traj} displays $\| \hat{R}(t) - I\|_F$ for both the S-TaSP and SL-EM solutions. The results show that the S-TaSP and SL-EM solutions are close to each other and both converge to $I$, which is consistent with the theoretical analysis. 

The results in Figure~\ref{fig:sdn_traj} demonstrate that the geometric deviation of $\hat R(t)$ is maintained at $O(10^{-13})$, illustrating the geometric preservation property of the S-TaSP and SL-EM algorithms.

Additionally, Figure~\ref{fig:time} plots the running time of the S-TaSP and SL-EM updates. The results show that the running time of S-TaSP is nearly $60\%$ that of SL-EM, highlighting its potential computational advantage in large-scale problems.

\section{Conclusion}

In this paper, we have developed a numerical scheme for simulating SDEs on the special orthogonal group $\SOn$. The proposed scheme is based on tangent space parametrization and avoids computing cost of the exponential map. We also provide a convergence analysis of the numerical scheme, and numerical simulation illustrates the effectiveness of our results.  

For future work, we aim to extend the stochastic tangent space parametrization scheme to general Riemannian manifolds and analyze the convergence performance.

\section{Appendix}
\subsection{Proof of Theorem~\ref{thm:correction}}
For the simplicity, we introduce some notations before the proof. For any $\lambda \in \R$, we denote $E(\lambda)$ as the $2\times2$ block $[\begin{smallmatrix} 0 & \lambda \\ -\lambda & 0 \end{smallmatrix}]$ and and $D(\lambda)$ as the $2\times2$ block $[\begin{smallmatrix} \lambda & 0 \\ 0 & \lambda \end{smallmatrix}]$. For any $\lambda_1,\lambda_2,\dots,\lambda_l$, and $n\ge 2l$, we denote the block-diagonal $\R^{n\times n}$ matrix as
\begin{align*}
    E_n(\lambda_1,\lambda_2,\dots,\lambda_l) &= \operatorname{diag}(E(\lambda_1),E(\lambda_2),\dots,E(\lambda_l),0), \\
    D_n(\lambda_1,\lambda_2,\dots,\lambda_l) &= \operatorname{diag}(D(\lambda_1),D(\lambda_2),\dots,D(\lambda_l),1).
\end{align*}

\noindent \textit{Proof:}
We write $ R = R_0 (I + \matZ + \bfC) $. Thus, it suffices to show that $ I + \matZ + \bfC \in \SOn $.

Since $ \matZ $ is skew-symmetric, by the real Schur decomposition, there exists an orthogonal matrix $ P $ and a sequence $\lambda_1,\lambda_2,\dots,\lambda_n \in \R$ such that:
\[
\matZ = P \Theta P^T, \quad \Theta = E_n(\lambda_1,\lambda_2,\dots,\lambda_l).
\]
Since $ \matZ^T \matZ = P \Theta^T \Theta P^T $, we have
\begin{align*}
   & I - \matZ^T \matZ = P (I - \Theta^T \Theta) P^T, \\ 
   & I - \Theta^T \Theta =D_n(1 - \lambda_1^2, 1 - \lambda_2^2, \dots,1-\lambda_l^2).
\end{align*}
Since $ I - \matZ^T \matZ $ is positive definite, we have $ 1 - \lambda_l^2 > 0 $ for all $ l $. Therefore, 
\begin{align*}
& \sqrt{I - \matZ^T \matZ} = P \Lambda P^T, \\ 
& \Lambda = D_n(\sqrt{1 - \lambda_1^2}, \sqrt{1 - \lambda_2^2},  \dots ,\sqrt{1 - \lambda_l^2}).
\end{align*}
We now examine $ I + \matZ + \bfC $. By expanding terms, we write
\[
I + \matZ + \bfC = P (\Theta + \Lambda) P^T.
\]
Since $ |\lambda_l| < 1 $, there exists an angle $ \theta_l \in [-\pi/2, \pi/2] $ such that $ \lambda_l = \sin(\theta_l) $ and $ \sqrt{1 - \lambda_l^2} = \cos(\theta_l) $. Setting
\[
\tilde{\matZ} = P \tilde{\Theta} P^T, \quad \tilde{\Theta} = E_n(\theta_1,\theta_2,\dots,\theta_l).
\]
it follows that $ \tilde{\matZ} \in \mathfrak{so}(n) $. By properties of the matrix exponential, we have
\[
\exp(\tilde{\matZ}) = P \exp(\tilde{\Theta}) P^T,
\]
where
\begin{align*}
     \exp(\tilde{\Theta}) = & D_n(\cos(\theta_1),\cos(\theta_2),\dots,\cos(\theta_l)) \\
    &+ E_n (\sin(\theta_1),\sin(\theta_2),\dots,\sin(\theta_l)) \\
    = & \Lambda + \Theta.
\end{align*}
Thus, $ \exp(\tilde{\matZ}) = I + \matZ + \bfC $, and since $ \tilde{\matZ} $ is skew-symmetric, we have $ \exp(\tilde{\matZ}) \in \SOn $, which completes the proof. \hfill $\blacksquare$
\subsection{Proof of Lemma~\ref{lem:1}}
We first introduce some lemmas.
\begin{lemma}
    Denote 
    \begin{align*}
        U_m = \sup_{ (m-1)\delta < t < m\delta} \| W(t) - W((m-1)\delta) \|_F
    \end{align*}
    is the maximum of the Frobenius norm of the Wiener process increment in the interval $[(m-1)\delta, m\delta]$. Then, there holds
    \begin{align*}
        \max_{m \in [M]} U^2_m \le \mathcal O(\delta^{1-\epsilon}),
    \end{align*}
    for abritrary small $\epsilon > 0$.
\end{lemma}
The following lemma gives the approximation behavior of the correction term $\bfC$.
\begin{lemma}
    For any $\matZ \in \mathfrak{so}(n)$ with $I - \matZ^T \matZ$ is positive definite, there holds
    \begin{align*}
       \| \bfC - \frac{1}{2} \matZ^2 \|_F \le \frac{1}{2} \| \matZ \|_F^4,
    \end{align*}
    where $\| \cdot \|_F$ is the Frobenius norm.
\end{lemma}
\noindent\textit{Proof:}
    Let $\matZ = P\Theta P^T$ be the Schur decomposition of $\matZ$ as in the proof of Theorem~\ref{thm:correction}. Then, we have
    \begin{align*}
    \begin{cases}
          \sqrt{I - \matZ^T \matZ} = P D_n(\sqrt{1 - \lambda_1^2},  \dots ,\sqrt{1 - \lambda_l^2}) P^T \\
         I + \frac{1}{2} \matZ^2 =  P D_n\left( 1 - \frac{1}{2} \lambda_1^2, \dots, 1 - \frac{1}{2} \lambda_l^2\right)  P^T \\
         \matZ^4 =   P D_n\left( \lambda_1^4, \dots, \lambda_l^4\right)  P^T 
    \end{cases}
    \end{align*}
    Thus, we have
    \begin{align*}
        \begin{cases}
            \bfC - \frac{1}{2} \matZ^2 =  P D_n\left( e_1, \dots, e_l\right)  P^T, \\
             e_i = \sqrt{1- \lambda_i^2} - 1 + \frac{1}{2}  \lambda_i^2.
        \end{cases}
    \end{align*}
  Since $I - \matZ^T \matZ$ is positive definite, we have $\theta_i < 1$ for all $i$. Thus, we have
    \begin{align*}
        -\frac{1}{2}  \lambda_i^4  \le  e_i = \sqrt{1- \lambda_i^2} - 1 + \frac{1}{2}  \lambda_i^2  \le 0,
    \end{align*}
    which implies that
    \begin{align*}
      \| \bfC  - \frac{1}{2} \matZ^2 \|_F =  & \sqrt{ 2  \sum_{i=1}^l \left( e_i \right)^2 } \le \sqrt{2 \sum_{i=1}^l (\frac{1}{2}  \lambda_i^4)^2} \\ 
      = & \frac{1}{2} \| \matZ^4 \|_F \le \frac{1}{2} \| \matZ \|_F^4,
    \end{align*}
    which completes the proof. \hfill $\blacksquare$\\
Here, we carry out the proof of Lemma~\ref{lem:1}. \\ 
\noindent{ \it Proof of Lemma~\ref{lem:1}.}  Denote $$B_{e,k} = B_{o,I}(\hat R_k) - \frac{1}{2} \sum_{j=1}^d B_j^2(\hat R_k).$$ First, we have
    \begin{align*}
        \hat R_m = & \hat R_{m-1} (I + \matZ_{m-1} + \bfC_{m-1}) \\
                 = & \hat R_{m-1} + \hat R_{m-1} (\matZ_{m-1} + \frac{1}{2} Z_{m-1}^2) +\hat  R_{m-1} E_{m-1} \\
                 = & \sum_{k=0}^{m-1} \hat R_k \left( \matZ_k + \frac{1}{2} \matZ_k^2 \right) + \sum_{k=0}^{m-1} \hat R_k E_k. \\
             E_k = & \bfC_k - \frac{1}{2} \matZ_k^2 = \mathcal O(\matZ_k^4).
    \end{align*}
    Thus, we have
    \begin{align*}
            & \hat R_m - R_{\delta}(m\delta) \\
          = & \sum_{k=0}^{m-1} \hat R_k \left( \matZ_k -  B_{e,k} \delta -  \sum_{j=1}^d B_j(\hat R_k) \sqrt{\delta} \epsilon_{j,k}  \right) \\ 
          & + \frac{1}{2} \sum_{k=0}^{m-1} \hat R_k \left( \matZ_k^2 - \sum_{j=1}^d B_j^2(\hat R_k) \right) + \sum_{k=0}^{m-1} \hat R_k E_k \\
          = & \hat R_k \sum_{k=0}^{m-1}( E^{b}_{k} + E^{\delta}_{k} ) + \sum_{k=0}^{m-1} \hat R_k E_k.
    \end{align*}
    where
    \begin{align*}
       & E^{b}_{k} = \frac{1}{2} \delta \hat R_k \sum_{k=0}^{m-1} \sum_{j = 1}^d B_j^2(\hat R_k) (\epsilon_{j,k}^2 - 1) \\ 
       & E^{\delta}_{k} = \frac{1}{2} \delta\hat R_k \sum_{k=0}^{m-1} \sum_{j' \neq j} B_j(\hat R_k) B_{j'}(\hat R_k) \epsilon_{j,k} \epsilon_{j',k}\\
       & + \frac{1}{2} \delta^{\frac{3}{2}} \hat R_k \sum_{k=0}^{m-1} B_{o,I}(\hat R_k) \sum_{j=1}^d B_j(\hat R_k) \epsilon_{j,k} + \frac{1}{2} \delta^2 \hat R_k B_{e,k}^2.
    \end{align*}
    Then we have, for $t \in [m\delta, (m+1)\delta]$,
    \begin{align*}
        & \hat R(t) - R_{\delta}(t) \\ 
        = & \hat R_m - R_{\delta}(m\delta) +  R_{\delta}(m\delta) - R_{\delta}(t) \\
        = & \hat R_k \sum_{k=0}^{m-1} E^{b}_{k}  + \sum_{k=0}^{m-1} \hat R_k E_k + (t - m\delta) \hat R_{m} B_{o,I}(\hat R_{m})  \\ 
        & + \hat R_k \sum_{k=0}^{m-1} E^{\delta}_{k} + \hat R_{m} \sum_{j=1}^d B_j(\hat R_{m}) \big( W_j(t)- W_j(m\delta) \big).
    \end{align*}
    Taking the expectation on both sides, we have
    \begin{align*}
        & \E \left[ \sup_{0 \le t \le T}\| \hat R(t) - R_{\delta}(t) \|_F^2 \right] \le 5a + 5b + 5c + 5d + 5e,    
    \end{align*}
    where
    \begin{align*}
        \begin{cases}
        a = \E \left[ \sup_{m \in [M]} \| \hat R_k \sum_{k=0}^{m-1} E^{b}_{k} \|_F^2 \right], \\
        b = \E \left[ \sup_{m \in [M]} \| \hat R_k \sum_{k=0}^{m-1} E^{\delta}_{k} \|_F^2 \right], \\
        c = \E \left[ \sup_{m \in [M]} \| \sum_{k=0}^{m-1} \hat R_k E_k \|_F^2 \right], \\
        d = \E \left[ \sup_{m,t \in \mathcal T} \| (t - m\delta) \hat R_{m} B_{o,I}(\hat R_{m}) \|_F^2 \right], \\
        e = \E \left[ \sup_{m,t \in \mathcal T} C^2 \| \hat R_{m} \sum_{j=1}^d \big( W_j(t)- W_j(m\delta) \big) \|_F^2 \right].
        \end{cases}
    \end{align*}
    where $\mathcal T = \{ m \in [M], t \in [m\delta, (m+1)\delta] \}$. 
    From~\cite{solo2021convergence}, we have $a, b, d = \calO(\delta), e = \calO(\delta^{1-\epsilon})$. The remaining part is to analyze $c$. From Doob's inequality, we have
    \begin{align*}
        &\E \left[ \sup_{m \in [M]} \| \sum_{k=0}^{m-1} \hat R_k E_k \|_F^2 \right]  \\ 
        \le & M \E \left[ \sum_{k=0}^{M-1} \| \hat R_k E_k \|_F^2 \right] 
        \le  9M \E \left[ \sum_{k=0}^{M-1} \| E_k \|_F^2 \right] 
    \end{align*}
    From~Lemma 2, we have
    \begin{align*}
        \| E_k \|_F^2 \le & \frac{1}{4} \| \matZ_k \|_F^8  
        \le \frac{1}{4} \left( \| B_{e,m}\delta + \sqrt{\delta} \sum_{j=1}^d B_j(\hat R_k) \epsilon_{j,k} \right)^8 \\
        \le & \frac{1}{4} \left( C \delta + \sqrt{\delta} C \sum_{j=1}^d | \epsilon_{j,k} | \right)^8 \\
                      \le & \frac{1}{4} \left( (d+1)^6 C^8 \delta^8 + (d+1)^6 C^8 \delta^4 \right) 
    \end{align*}
    which implies that
    \begin{align*}
        \E \left[ \sup_{m \in [M]} \| \sum_{k=0}^{m-1} \hat R_k E_k \|_F^2 \right]  \le M^2 O(\delta^4) = o(\delta).
    \end{align*}
    
    Finally, we have
    \begin{align*}
        &\E \left[ \sup_{0 \le t \le T}\| \hat R(t) - R_{\delta}(t) \|_F^2 \right] \\ 
        \le & 5a + 5b + 5c + 5d + 5e = \calO(\delta^{1-\epsilon}),
    \end{align*}
    which completes the proof. \hfill $\blacksquare$

\bibliographystyle{IEEEtran}
\bibliography{ref}

\end{document}